\documentclass{amsproc}
\usepackage{amscd, amsfonts, amsmath,amssymb,amsthm, mathtools,mathrsfs}
\usepackage[normalem]{ulem}
\usepackage{graphicx}
\usepackage[cmtip,all]{xy}
\copyrightinfo{2009}{American Mathematical Society}

\newtheorem{theorem}{Theorem}[section]
\newtheorem{corollary}[theorem]{Corollary}
\newtheorem{proposition}[theorem]{Proposition}

\theoremstyle{definition}
\newtheorem{definition}[theorem]{Definition}
\newtheorem{example}[theorem]{Example}

\theoremstyle{remark}
\newtheorem{remark}[theorem]{Remark}

\numberwithin{equation}{section}

\newcommand{\kk}{\Bbbk}

\newcommand{\Z}{\mathbb{Z}}

\newcommand{\N}{\mathbb{N}}

\newcommand{\HH}{\mathrm{HH}}

\DeclareMathOperator{\GL}{GL}

\title[Twists in algebra, representation theory, and geometry]{A Primer on Twists in the Noncommutative Realm Focusing on Algebra, Representation Theory, and Geometry}

\author[Ocal]{Pablo S.\ Ocal}
\address{(Ocal) UCLA Mathematics Department, Los Angeles, CA 90095, U.S.A.}
\email{socal@math.ucla.edu}
\thanks{}

\author[Ueyama]{Kenta Ueyama}
\address{(Ueyama) Department of Mathematics, Faculty of Science, Shinshu University, Matsumoto, Nagano 390-8621, Japan}
\email{ueyama@shinshu-u.ac.jp}

\author[Veerapen]{Padmini Veerapen}
\address{(Veerapen) Department of Mathematics, Tennessee Tech University, Cookeville, TN 38505, U.S.A.}
\email{pveerapen@tntech.edu}



\begin{document}

\begin{abstract}
We review several techniques that twist an algebra's multiplicative structure. We first consider twists by an automorphism, also known as Zhang twists, and we relate them to 2-cocycle twists of certain bialgebras. We then outline the classification and properties of twisted tensor products, and we examine twisted Segre products. Our exposition emphasizes clarity over generality, providing a wealth of interconnecting examples.
\end{abstract}

\maketitle

\section{Introduction}
Twisting the multiplicative structure of algebras and of coalgebras has long been an area of interest in mathematics; some classic examples include skew polynomial rings, skew group rings, and Ore extensions. Here we present recent constructions that generalize those.

In the early 1990s, the notion of a twist of a graded algebra by an automorphism was introduced by Artin, Tate, and Van den Bergh in \cite{ATV1991}, and was later generalized by Zhang in \cite{Zhang1996}. This twist of a graded algebra became known as a Zhang twist, and it has played a pivotal role in various areas of noncommutative algebra and noncommutative projective geometry. Many fundamental properties of graded algebras are preserved under a Zhang twist, such as their graded module categories, noncommutative projective schemes, Gelfand-Kirillov dimension, global dimension, graded $\operatorname{Ext}$ groups, and Artin-Schelter regularity \cite{Ro2016,Zhang1996}.

At the same time, similar notions were being developed for bialgebras and Hopf algebras. In the late 1980s, motivated by the tensor equivalence of two module categories, Drinfel'd \cite{Dr87} introduced a twist for the coalgebra structure of a Hopf algebra. These are known as Drinfel'd twists, and similarly to Zhang twists, they preserve important properties of Hopf algebras \cite{AEGN}. For instance, certain classes of cosemisimple unimodular, cosemisimple involutive, and cosemisimple quasitriangular finite dimensional Hopf algebras are invariant under Drinfel'd twists. The dual version of the Drinfel'd twist, called a 2-cocycle twist, was studied by Doi and Takeuchi \cite{Doi93,DT94} in the early 1990s. These modify the algebra structure of a Hopf algebra and yield a tensor equivalence of the corresponding comodule categories. We focus on Zhang twists and 2-cocycle twists of $\mathbb{Z}$-graded algebras in Section \ref{SecZhangtwist}, but there are other notions of twists of noncommutative algebras (see e.g., \cite{Davies2017, M2005}).

Extending further the idea of modifying existing multiplicative structures, we consider twisted tensor products of algebras in Section \ref{Sec3}. This notion was introduced by \v{C}ap, Schichl, and Van\v{z}ura \cite{MR1352565}, it is contemporary to Zhang twists and Drinfel'd twists, and when graded algebras are twisted, it inherits a bigrading. Twisted tensor products were originally designed to provide an appropriate algebraic interpretation of the noncommutative product of topological spaces, but since they encompass vast families of algebras, they have also proven useful to study representation theory. For example, Ore extensions, noncommutative $2$-tori, crossed products of $C^*$-algebras with groups, and universal enveloping algebras of Lie algebras, are all particular cases of twisted tensor products. Because of their generality, it is difficult to discern which properties or structures are preserved by twisted tensor products, and little is known besides Artin-Schelter regularity and Frobenius structures \cite{MR3866681, OcalOswald}.

As an application of twisted tensor products, in Section \ref{Sec4} we consider twisted Segre products, introduced by the second author and collaborator in \cite{HU}. These are a noncommutative generalization of the Segre product in commutative algebra and algebraic geometry. In particular, we present a categorical relationship between twisted tensor products and twisted Segre products.

\subsection*{Convention} Throughout the paper, let $\kk$ be a field. An \emph{algebra} means a unital associative $\kk$-algebra and unadorned $\otimes$ means $\otimes_\kk$. For an algebra $A$ we denote its unit map by $\eta_A: k \to A$ and its multiplication map by $\nabla_A: A \otimes A \to A$. For a coalgebra $A$ we denote its counit map by $\epsilon_A: A \to k$ and its comultiplication map by $\Delta_A: A \to A \otimes A$. When $A$ is graded, we denote by $|a|$ the degree of a homogeneous element $a\in A$. We denote the composition of maps by their concatenation.

\section{Zhang Twists of Algebras and 2-cocycle Twists of Bialgebras}\label{SecZhangtwist}
The notion of a twist of a $\mathbb Z$-graded algebra was first introduced in \cite[Section 8]{ATV1991}. Using the more general notion of twisting systems, Zhang proved in \cite{Zhang1996} that given two $\mathbb N$-graded algebras $A$ and $B$ generated in degree one, $A$ is isomorphic to a twisted algebra of $B$ if and only if the graded module categories of $A$ and $B$ are equivalent. 

\begin{definition} 
\label{defn:Zhang-twist} 
Let $A$ be a $\mathbb{Z}$-graded algebra and $\phi$ be a graded degree-zero automorphism of $A$. The \emph{right Zhang twist} $A^{\phi}$ of $A$ by $\phi$ is the algebra that coincides with $A$ as a graded vector space over $\kk$, with the multiplication
\[r *_\phi s \coloneqq r \phi^{|r|}(s),
\quad \text{for any homogeneous elements } r,s \in A. \]
The \emph{left Zhang twist} $\prescript{\phi}{}{A}$ of $A$ by $\phi$ is defined similarly, with the multiplication
\[r *_\phi s \coloneqq \phi^{|s|}(r) s,
\quad \text{for any homogeneous elements } r,s \in A. \]
\end{definition} 

For $A$ and $\phi$ given in Definition \ref{defn:Zhang-twist}, the left and right Zhang twist of $A$ by $\phi$ are again graded associative algebras \cite[Proposition 4.2]{Zhang1996}. Furthermore, there is a natural isomorphism $A^\phi \cong \prescript{\phi^{-1}}{}{A}$ mapping $a$ to $\phi^{-|a|}(a)$ \cite[Lemma 4.1]{Lecoutre-Sierra2019}.

\begin{example}\label{exam.zhang-a}
Consider the polynomial ring $\kk[x_1,x_2,x_3,x_4]$ and define a graded degree-zero automorphism $\theta:\kk[x_1,x_2,x_3,x_4] \to \kk[x_1,x_2,x_3,x_4]$ by
\begin{equation*}
\theta(x_1) = ax_1,\quad \theta(x_2) = ax_2, \quad \theta(x_3) = x_3,\quad \theta(x_4) = x_4
\end{equation*}
for some $a\in\kk^{\times}$. The Zhang twist $\kk[x_1,x_2,x_3,x_4]^{\theta}$ can be seen as the quotient of the free algebra on generators $x_1, x_2, x_3$, and $x_4$ subject to the relations
\begin{align*}
&x_1x_2 - x_2x_1,  && ax_1x_3 - x_3x_1, && ax_1x_4 - x_4x_1,\\
&ax_2x_3 - x_3x_2, && ax_2x_4 - x_4x_2, && x_3x_4 - x_4x_3.
\end{align*}
This will be a running example throughout the paper.
\end{example}
Zhang twists may also be defined in the more general setting of an algebra graded by a semigroup by considering twisting systems \cite{Zhang1996}. 

\begin{definition}  \cite[Definition 2.1]{Zhang1996}
Let $A = \bigoplus_{n \in \Z}A_n$ be a $\Z$-graded algebra and let $\tau = \{\tau_n | n \in \Z\}$ be a set of $\kk$-linear automorphisms of $A$. We call $\tau$ a \emph{twisting system} when
\[\tau_n(r\tau_m(s)) = \tau_{n}(r)\tau_{n+m}(s)\]
for all $n, m, \ell \in \Z$ and $r \in A_m, s \in A_{\ell}$. Given a twisting system, we may define a new graded multiplication $\ast$ on the underlying graded $\kk$-module $\bigoplus_{n \in \Z}A_n$ by 
\[r \ast s = r \tau_m(s)\] for all $r \in A_m$, $s \in A_{\ell}$.
\end{definition} 

In the case of a $\mathbb{Z}$-graded algebra, such a twisting system arises from a single algebra automorphism of the associated $\mathbb{Z}$-algebra \cite[Proposition 4.2]{S2011}.

\begin{example}\cite{Tran}
Suppose char$(\kk) \ne 2$. Consider the following family of quadratic algebras
\[A(\rho) = \kk \langle x_1, \dots, x_4 \rangle / (w_1, \dots, w_6), \]
where $\rho \in \kk^{\times}$ and
\begin{align*}
w_1 &= \rho^2x_1x_2+x_2x_1, & w_2 &= \rho x_1x_3 - x_3x_1, & w_3 &= \rho x_1x_4 + x_4x_1,\\
w_4 &= x_2x_3 - \rho x_3x_2, & w_5 &= x_2x_4 + \rho x_4 x_2, & w_6 &= x_3x_4 + x_4x_3 + x_1x_2.
\end{align*}
This is a family of quadratic Artin-Schelter regular algebras. It is noteworthy that the algebras in this family are twists of each other via a twisting system. In particular, the twisting system cannot be reduced to a twist by a single algebra automorphism. 
\end{example}

We now consider Zhang twists on algebras and bialgebras (respectively Hopf algebras), where one should note that Zhang twists only deform the algebra structure. With this in mind, we let $B$ be a $\mathbb Z$-graded bialgebra (respectively Hopf algebra). We say that an endomorphism $\phi\in {\rm End}_\kk(B)$ is a \emph{graded bialgebra automorphism} (respectively ~\emph{graded Hopf algebra automorphism}) of $B$ if it is a bialgebra (respectively Hopf algebra) automorphism of $B$ preserving the grading of $B$, that is $\phi(B_n)\subseteq B_n$ for all $n\in \mathbb Z$.

The rest of the section is devoted to relating Zhang twists to 2-cocycle twists using recent work in \cite{HNUVVW} by the third author and collaborators. First, we review 2-cocycle twists of bialgebras and of Hopf algebras. As mentioned in the introduction, 2-cocycle twists originated from the work of Doi and Takeuchi \cite{Doi93, DT94} as a dual version to the twist defined by Drinfel'd \cite{Dr87}.

\begin{definition}\cite[Section 1]{M2005}
\label{M2005}
Suppose $H$ is a Hopf algebra. Let $\sigma : H \otimes H \to \kk$ be a 2-cocycle for $H$. The bialgebra $H^{\sigma}$ is called a \emph{cocycle twist} of $H$ with respect to $\sigma$ when
\begin{enumerate}
\item $H^{\sigma} = H$ as a coalgebra,
\item $H^{\sigma}$ has multiplication given by
\[h \cdot_{\sigma} l \coloneqq \sum \sigma^{-1}(h_1, l_1)h_2l_2\sigma(h_3, l_3).\]
\end{enumerate}
The bialgebra $H^{\sigma}$ admits a Hopf algebra structure with the antipode
\[S^{\sigma}(h) \coloneqq \sum \sigma^{-1}(h_1, Sh_2)Sh_3\sigma(Sh_4, h_5).\]
\end{definition}

We note that the notation $H^{\sigma}$ is a modification of the one used in \cite{M2005}. Now, to relate a Zhang twist to a 2-cocycle twist, we introduce the notion of \textit{twisting pairs}.

\begin{definition}\cite[Definition D]{HNUVVW}
Let $B$ be a bialgebra. A pair $(\phi_1, \phi_2)$ of algebra automorphisms of $B$ is said to be a \emph{twisting pair} when the following conditions hold.
\begin{enumerate}
\item $\Delta \phi_1 = (1 \otimes \phi_1) \Delta$ and $\Delta \phi_2 = (\phi_2 \otimes  1) \Delta$, \mbox{ and }
\item $\varepsilon \phi_1 \phi_2 = \varepsilon$.
\end{enumerate}
\end{definition}

One may lift a bialgebra $B$ to a Hopf algebra $\mathcal{H}(B)$ using Takeuchi's Hopf envelope construction \cite{T71}. Moreover, if $B$ is $\Z$-graded as an algebra and $\Delta(B_n) \subset B_n \otimes B_n$ for all $n \in \Z$, then $\mathcal{H}(B)$ also satisfies these conditions. We use Takeuchi's explicit construction of $\mathcal{H}(B)$ in the result below.

\begin{theorem}\cite[Theorem 2.3.3]{HNUVVW}
\label{Sec2BigThm}
Let $B$ be a bialgebra that is $\Z$-graded as an algebra and satisfying $\Delta(B_n) \subset B_n \otimes B_n$ for all $n \in \Z$. Given a twisting pair $(\phi_1, \phi_2)$ of $B$ there is a unique pair $(\mathcal{H}(\phi_1), \mathcal{H}(\phi_2))$ of the Hopf envelope of $\mathcal{H}(B)$ extending $(\phi_1, \phi_2)$. Moreover, the 2-cocycle twist $\mathcal{H}(B)^{\sigma}$, with 2-cocycle $\sigma: \mathcal{H}(B) \otimes \mathcal{H}(B) \to \kk$ given by 
\[\sigma(x, y) = \varepsilon(x)\varepsilon(\mathcal{H}(\phi_2)^{|x|}(y))\quad \mbox{for any homogeneous elements } x, y \in \mathcal{H}(B),\]
is the right Zhang twist $\mathcal{H}(B)^{\mathcal{H}(\phi_1)\mathcal{H}(\phi_2)}$.
\end{theorem}

\vspace{1.5mm}

We now provide an application of the above construction by computing twisting pairs explicitly for $\mathcal O_q(M_n(\kk))$, the one parameter quantization of the coordinate ring of $n \times n$ matrices.

\begin{example}\label{ex:gl}
Let $n\ge 2$ be an integer and let $q\in \kk^\times $ be any scalar. As an algebra, $\mathcal O_q(M_n(\kk))$ is generated by $n^2$-generators $\{x_{ij}\}_{1\leq i,j\leq n}$ subject to the relations
 \begin{equation*}
     \label{eq:Rqrelation1}
     \left\{\begin{aligned}
qx_{ks} x_{us} &= x_{us} x_{ks}&&{\rm if }\,\, k<u\\ 
qx_{ks} x_{kv} &= x_{kv} x_{ks} &&{\rm if }\,\,  s<v\\
x_{us} x_{kv} &= x_{kv} x_{us}&&{\rm if }\,\,  s<v, k<u\\
x_{us} x_{kv} &= x_{kv} x_{us}+(q-q^{-1})x_{ks} x_{uv} &&{\rm if }\,\, s<v,
u<k.
 \end{aligned}\right.
 \end{equation*}
 The coalgebra structure on $\mathcal O_q(M_n(\kk))$ given by
 \begin{equation*}
     \Delta(x_{ij})=\sum_{1\leq k\leq n} x_{ik}\otimes x_{kj}\quad \text{and}\quad \varepsilon(x_{ij})=\delta_{ij}\quad \text{for all}\ 1\leq i,j\leq n,
 \end{equation*}
 makes $\mathcal O_q(M_n(\kk))$ a bialgebra. Our goal is to lift $\mathcal O_q(M_n(\kk))$ to its Hopf envelope $\mathcal O_q(\GL_n(\kk))$, the one-parameter quantization of the coordinate ring of the general linear group. Following \cite{T2002}, this is achieved by first defining a coquasitriangular structure on $\mathcal O_q(M_n(\kk))$ and then localizing at a central grouplike element. As it were, $\mathcal{O}_q(M_n(\kk))$ is coquasitriangular as it is isomorphic to the coquasitriangular bialgebra $A(R_q)$ obtained via the Faddeev-Reshetikhin-Takhtajan construction (see e.g., \cite{FRT,CMZ}). For any scalar $q \in \kk^{\times}$, the classical Yang-Baxter operator $R_q$ on $V$ (cf. \cite{T2002}) is given by  
\begin{equation*}
    R_q(v_i\otimes v_j)=\begin{cases}
    qv_i\otimes v_i &\qquad i=j\\
     v_i\otimes v_j &\qquad i<j\\
    v_i\otimes v_j+(q-q^{-1})v_j \otimes v_i &\qquad i>j,
    \end{cases}
\end{equation*}
where $V$ is an $n$-dimensional vector space with basis $\{v_1,\ldots,v_n\}$. The central group-like element $g$ of $\mathcal O_q(M_n(\kk))$ given by
\begin{equation*}
\label{central g}
  g\coloneqq \sum_{\sigma\in S_n}(-q)^{-l(\sigma)}x_{\sigma(1)1} \cdots x_{\sigma(n)n}
\end{equation*}
where $l(\sigma)$ denotes the length of the permutation $\sigma$, is called the $q$-determinant.
\end{example}

\begin{proposition}\cite[Proposition 3.1.7]{HNUVVW}
All twisting pairs of $\mathcal O_q(\GL_n(\kk))$ are of the form $(\phi_1, \phi_2)$, where $\phi_1$ and $\phi_2$ are given by
 \[
 \begin{pmatrix*}[c]
\phi_1(x_{11})         & \cdots  & \phi_1(x_{1n}) \\
\vdots               & \ddots   & \vdots             \\
\phi_1(x_{n1})         & \cdots  & \phi_1(x_{nn}) 
\end{pmatrix*}=
 \begin{pmatrix*}[c]
x_{11}         & \cdots  & x_{1n} \\
\vdots               & \ddots   & \vdots             \\
x_{n1}         & \cdots  & x_{nn} 
\end{pmatrix*}
\begin{pmatrix*}[c]
\alpha_{11}         & \cdots  & \alpha_{1n} \\
\vdots               & \ddots   & \vdots             \\
\alpha_{n1}         & \cdots  & \alpha_{nn} 
\end{pmatrix*}
 \]
 and 
 \[
 \begin{pmatrix*}[c]
\phi_2(x_{11})         & \cdots  & \phi_2(x_{1n}) \\
\vdots               & \ddots   & \vdots             \\
\phi_2(x_{n1})         & \cdots  & \phi_2(x_{nn}) 
\end{pmatrix*}= \begin{pmatrix*}[c]
\alpha_{11}         & \cdots  & \alpha_{1n} \\
\vdots               & \ddots   & \vdots             \\
\alpha_{n1}         & \cdots  & \alpha_{nn} 
\end{pmatrix*}^{-1}
 \begin{pmatrix*}[c]
x_{11}         & \cdots  & x_{1n} \\
\vdots               & \ddots   & \vdots             \\
x_{n1}         & \cdots  & x_{nn} 
\end{pmatrix*},
 \]
with $(\alpha_{ij})\in \GL_n(\kk)$. Moreover, we have the following.
\begin{enumerate}
\item If $q=1$, then the maps $(\phi_1, \phi_2)$ defined above form a twisting pair for all $(\alpha_{ij}) \in \GL_n(\kk)$.
\item If ${\rm char}(\kk)\neq 2$ and $q=-1$, then $(\alpha_{ij})$ defines a twisting pair as above if and only if it is a generalized permutation matrix. 
\item If $q \neq \pm 1$, then $(\alpha_{ij})$ defines a twisting pair as above if and only if it is is diagonal.
\end{enumerate}
 In particular, $\phi_1(g^{-1})=(-1)^{l(\tau)}|(\alpha_{ij})|^{-1}g^{-1}$, $\phi_2(g^{-1})=(-1)^{l(\tau)}|(\alpha_{ij})|g^{-1}$, and $\phi_1 \phi_2(g^{-1})=g^{-1}$, where $\tau$ is trivial if $q\ne -1$, and $\tau\in S_n$ satisfied $|(\alpha_{ij})|=(-1)^{l(\tau)}\alpha_{\tau(1)1}\cdots \alpha_{\tau(n)n}$ if $q=-1$. 
\end{proposition}

The above twisting pairs of $\mathcal O_q(\GL_n(\kk))$ allow us to apply Theorem \ref{Sec2BigThm}.

\begin{corollary}
\cite{HNUVVW}
The right Zhang twist of $\mathcal O_q(\GL_n(\kk))$ by $\phi_1\phi_2$ is isomorphic to the 2-cocycle twist of $\mathcal O_q(\GL_n(\kk))$ by the 2-cocycle $\sigma$ satisfying
\begin{align*}  
\sigma &(x_{i_1j_1}^{p_1}\cdots x_{i_sj_s}^{p_s}g^{-r}, x_{u_1v_1}^{q_1}\cdots x_{u_kv_k}^{q_k}g^{-t})\\
&= (-1)^{mtl(\tau)}\delta_{i_1, j_1}\cdots \delta_{i_s, j_s}\left(((\alpha_{ij})^{-m})_{u_1v_1}\right)^{q_1} \cdots \left(((\alpha_{ij})^{-m})_{u_kv_k}\right)^{q_k}|(\alpha_{ij})|^{mt},
\end{align*}
where $m=p_1+\dots +p_s-nr$ and $1\leq i_\bullet,j_\bullet,u_\bullet,v_\bullet\leq n$ and $p_\bullet, q_\bullet,s,k, r,t$ are non-negative integers, $\tau$ is trivial if $q\ne -1$, and $\tau\in S_n$ satisfies $|(\alpha_{ij})|=(-1)^{l(\tau)}\alpha_{\tau(1)1}\cdots \alpha_{\tau(n)n}$ if $q=-1$.
\end{corollary}

\section{Twisted Tensor Products}\label{Sec3}

When twisting the multiplication of a given algebra, a natural option is to consider noncommutative analogues of tensor products. Such an analogue was introduced by \v{C}ap, Schichl, and Van\v{z}ura \cite{MR1352565}, and is known as a \emph{twisted tensor product}. These are ubiquitous in noncommutative algebra and noncommutative geometry.

\begin{definition} \cite[Definition 2.3]{MR1352565}
Let $A$ and $B$ be algebras. A \emph{twisting map} $\tau:B\otimes A\rightarrow A\otimes B$ is a bijective $\kk$-linear map satisfying
\begin{align*}
1\otimes \eta_B &= \tau (\eta_B\otimes 1),\\
\eta_A\otimes 1 &= \tau (1\otimes \eta_A),\\
\tau (\nabla_B \otimes \nabla_A) &= (\nabla_A \otimes
\nabla_B)(1\otimes \tau \otimes 1)(\tau \otimes \tau) (1\otimes \tau \otimes 1).
\end{align*}
\end{definition}

The first two conditions express the compatibility of $\tau$ with units $\eta_A$ and $\eta_B$, while the third ensures the compatibility of $\tau$ with the multiplication maps $\nabla_A$ and $\nabla_B$. The commutative diagrams given by the above equalities are a convenient computational tool. Some authors allow twisting maps to be non-invertible, but that will not be our case. When working with graded algebras $A = \bigoplus_{i\in\Z} A_i$ and $B = \bigoplus_{i\in\Z} B_i$, it is often desirable to require that the twisting map be \emph{strongly graded}, namely $\tau(B_j \otimes A_i) \subseteq A_i \otimes B_j$.

\begin{definition} \cite[Definition 2.1]{MR1352565}
Let $A$ and $B$ be algebras and let $\tau:B\otimes A\rightarrow A\otimes B$ be a twisting map. The \emph{twisted tensor product} $A\otimes_{\tau} B$ is the $\kk$-vector space $A\otimes B$ with the following unit and multiplication maps.
\begin{align*}
\eta_{A\otimes_{\tau} B} &\coloneqq (\eta_A\otimes \eta_B),\\
\nabla_{A\otimes_{\tau} B} &\coloneqq (\nabla_A\otimes \nabla_B)(1\otimes \tau \otimes 1).
\end{align*}
\end{definition}

\begin{proposition} \cite[Proposition 2.3]{MR1352565}
Let $A$ and $B$ be algebras and let $\tau:B\otimes A\rightarrow A\otimes B$ be a twisting map. Then $(A\otimes_{\tau} B,\nabla_{A\otimes_{\tau} B},\eta_{A\otimes_{\tau} B})$ is an algebra.
\end{proposition}

In other words, the compatibility of $\tau$ with the unital associative algebra structures of $A$ and $B$ indeed guarantees that $A\otimes_{\tau} B$ inherits this structure. The compatibility of $\tau$ with $\nabla_A$ and $\nabla_B$ can be interpreted as the equalities
\begin{align*}
\tau(1\otimes \nabla_A) &= (\nabla_A\otimes 1)(1\otimes \tau)(\tau\otimes 1),\\
\tau(\nabla_B\otimes 1) &= (1\otimes \nabla_B)(\tau\otimes 1)(1\otimes \tau).
\end{align*}
Intuitively, the above equalities indicate that $\tau$ commutes with the multiplication maps. It turns out that an algebra $\Lambda$ can be expressed as a twisted tensor product of two subalgebras $A$ and $B$ essentially when $\Lambda$ is isomorphic to $A\otimes B$ as vector spaces, and surprisingly the twisting map is unique.

\begin{theorem} \cite[Proposition 2.7]{MR1352565}
Let $\Lambda$, $A$, $B$ be algebras. Then $\Lambda \cong A\otimes_{\tau} B$ for some twisting map $\tau : B\otimes A\to A\otimes B$ if and only if there are injective algebra morphisms $\iota_A : A\to \Lambda$ and $\iota_B : B\to \Lambda$ such that $\nabla_{\Lambda} (\iota_A\otimes\iota_B) : A\otimes B\to \Lambda$ is an isomorphism of $\kk$-vector spaces. Moreover
\begin{equation*}
\tau = (\nabla_{\Lambda} (\iota_A\otimes\iota_B))^{-1} \nabla_{\Lambda} (\iota_B\otimes\iota_A)
\end{equation*}
whence $\tau$ is unique.
\end{theorem}

Note that $\iota_A : A\to \Lambda$ and $\iota_B : B\to \Lambda$ are then given by $\iota_A(a) = a\otimes 1_B$ and $\iota_B(b) = 1_A\otimes b$. Moreover, twisted tensor products inherit a universal property that recovers the usual universal property of the usual tensor product when the twisting map is trivial.

\begin{theorem}
Let $A$, $B$ be algebras, let $\tau : B\otimes A\to A\otimes B$ be a twisting map, let $\iota_A : A \to A\otimes_{\tau} B$ and $\iota_B : B \to A\otimes_{\tau} B$ be such that $\tau = \nabla_{A\otimes_{\tau} B} (\iota_A \otimes \iota_B)$. If $\Lambda$ is an algebra and $i_A : A\to \Lambda$ and $i_B : B \to \Lambda$ are algebra morphisms such that $\nabla_{\Lambda} (i_B\otimes i_A) = \nabla_{\Lambda} (i_A\otimes i_B) \tau$ then there is a unique algebra morphism $h:A\otimes_{\tau} B\to \Lambda$ satisfying $h \iota_A = i_A$ and $h \iota_B = i_B$.
\end{theorem}

In other words, the triple $(A,A\otimes_{\tau} B,B)$ can be seen as an initial object in the class of cospans $A \to \Lambda \leftarrow B$ of algebras having varying apex $\Lambda$, fixed feet $A$ and $B$, and commute with $\tau$ in the appropriate sense. As with the usual tensor product, twisted tensor products can be iterated, and the iterations will inherit the aforementioned universal property.

\begin{theorem} \cite[Theorem 2.1]{MR2458561}
Let $A$, $B$, $C$ be algebras and let $\tau_{AB} : B\otimes A\to A\otimes B$, $\tau_{BC} : C\otimes B\to B\otimes C$, $\tau_{AC} : C\otimes A\to A\otimes C$ be twisting maps. Then $(1\otimes \tau_{BC})(\tau_{AC}\otimes 1) : C\otimes (A\otimes_{\tau_{AB}} B)\to (A\otimes_{\tau_{AB}} B)\otimes C$ is a twisting map if and only if $(\tau_{AB}\otimes 1)(1\otimes \tau_{AC}) : (B\otimes_{\tau_{BC}} C)\otimes A\to A\otimes (B\otimes_{\tau_{BC}} C)$ is a twisting map. In that case we have the isomorphism of algebras
\begin{equation*}
(A\otimes_{\tau_{AB}} B)\otimes_{(1\otimes \tau_{BC})(\tau_{AC}\otimes 1)} C \cong A\otimes_{(\tau_{AB}\otimes 1)(1\otimes \tau_{AC})} (B\otimes_{\tau_{BC}} C).
\end{equation*}
\end{theorem}

This gives one additional compatibility condition on the twisting maps, but as long as it is satisfied, the twisted tensor product of a family of algebras can be found by consecutive forming twisted tensor products.

\begin{theorem} \cite[Theorem 2.9]{MR2458561}
Fix $n \in\N$. Let $A_i$ be algebras for all $i = 1,\dots, n$, and let $\tau_{ij} : A_j\otimes A_i\to A_j\otimes A_i$ be twisting maps for all $1\leq i<j\leq n$. If $\tau_{ij}$, $\tau_{jk}$, and $\tau_{ik}$ are compatible for all $1\leq i<j<k\leq n$ then the iterated twisted tensor product $A_1\otimes_{\tau_{12}} \cdots \otimes_{\tau_{i-1 i}} A_i\otimes_{\tau_{i i+1}} \cdots \otimes_{\tau_{n-1 n}} A_n$ is well defined.
\end{theorem}

Twisted tensor products present a higher level of complexity when compared to their commutative counterparts. Of course, when the twisting map transposes the elements, namely $\tau:B\otimes A\to A\otimes B$ is given by $\tau(b\otimes a) = a\otimes b$ for all $a\in A$ and $b\in B$, we recover the usual tensor product, but classifying the other outcomes for most twisting maps is a gargantuan task. As a result of its universal property, a twisted tensor product is uniquely determined (whence in bijection) by a twisting map, whence it is enough to classify twisting maps \cite[Theorem 2.4]{MR2235816}. The first approaches to a classification consisted essentially of treating twisted tensor products on a case-by-case basis. For example, one can see that the available options for $A\otimes_{\tau} \kk[x]$ coincide precisely with the Ore extensions of $A$.

\begin{example}
Consider an algebra $A$, an algebra automorphism $\sigma : A\to A$, and a $\kk$-linear map $\delta:A\to A$  such that $\delta(\nabla(r\otimes s)) = \nabla(\delta(r)\otimes s) + \nabla(\sigma(r)\otimes \delta(s))$ for all $r,s\in A$. We define a twisting map $\tau:\kk[x]\otimes R \to R\otimes \kk[x]$ by setting $\tau(x\otimes r) = \sigma(r)\otimes x + \delta(r)\otimes 1$, and extend by linearity, together with imposing the compatibility conditions to compute $\tau(x^i\otimes r)$ for all $i\in\N$. The resulting twisted tensor product is isomorphic to the Ore extension $A[x;\sigma,\delta]$, namely the vector space $A[x]$ with multiplication $xr = \sigma(r)x + \delta(r)$ for all $r\in A$ (see \cite{MR1321145}).
\end{example}

The case $A\otimes_{\tau} \kk^2$ has a similar flavor to those of Ore extensions.

\begin{theorem} \cite[Proposition 2.10]{MR2235816}
Let $A$ be an algebra, let $\sigma : A \to A$ be an algebra automorphism, let $\prescript{}{\sigma}{A}$ be the $A$ bimodule with underlying vector space $A$ and action $\rho : A\otimes \prescript{}{\sigma}{A} \otimes A \to \prescript{}{\sigma}{A}$ given by $\rho \coloneqq \nabla (\sigma \otimes \nabla)$, let $\delta : A \to \prescript{}{\sigma}{A}$ be an idempotent derivation satisfying $\delta = \delta^2 + \sigma \delta + \delta \sigma$. The twisted tensor products $A\otimes_{\tau} \kk^2$ are in bijection with the pairs $(\sigma,\delta)$.
\end{theorem}

There are only partial results for the case of algebras with truncated polynomials \cite{MR2927170}. For $\kk[y]\otimes_{\tau}\kk[x]/(x^2)$ the classification is already considerably harder.

\begin{theorem} \cite[Section 4]{MR2642035}
Let $\kk$ be a field of characteristic zero. A choice of $m\in 2\N$, $p_0,q_m\in \kk^{\times}$, and $q_0,\dots,q_{m-1}\in \kk$ uniquely determines the polynomials $p(y) = -y + p_0$ and $q(y) = q_0 + q_1 y + \cdots + q_m y^m$, corresponding to a unique twisting map $\tau: \kk[x]/(x^2)\otimes_{\tau} \kk[y] \to \kk[y]\otimes_{\tau}\kk[x]/(x^2)$ given as $\tau(x\otimes a(y)) = \iota_0(a(y))\otimes 1 + \iota_1(a(y))\otimes x$ where $\iota_0(y) = q(y)$, $\iota_1(y) = p(y)$, $\iota_1(a(y)b(y)) = \iota_1(a(y))\iota_1(b(y))$, $\iota_0(a(y)b(y)) = \iota_0(a(y))b(y) + \iota_1(a(y))\iota_0(b(y))$, $\iota_0(a(y))^2 = 0$, and $\iota_0 \iota_1(a(y)) = -\iota_1 \iota_0(a(y))$ for all $a(y), b(y)\in \kk[y]$.
\end{theorem}

Some partial results have been achieved for $\kk^n$ with $\kk^m$ where $n,m\in\N\setminus \{0\}$ \cite{MR4034796, MR2832263}. These rely on rewriting the twisting map $\tau:\kk^m \otimes \kk^n\to \kk^n\otimes \kk^m$ as a family of matrices $\mathcal{M}_{\tau}$, translating the classification to a linear algebra problem. Given $e_1,\dots,e_m$ and $f_1,\dots,f_n$ the canonical basis of $\kk^m$ and $\kk^n$ respectively, we can write $\tau(e_i\otimes f_j) = \sum_{r,s}{\lambda_{ij}^{rs} f_r\otimes e_s}$. Set $M_{\tau}(i,s)_{rs} = \lambda_{ij}^{rs}$ and $\mathcal{M}_{\tau} = \{M_{\tau}(i,s)\}_{r,j}$ for $i,s = 1,\dots,m$, $j,r = 1,\dots,n$. Then $\tau$ is a twisting map if and only if $M_{\tau}(1,s),\dots,M_{\tau}(m,s)$ are orthogonal idempotent for all $s = 1,\dots,m$, and the dual condition exchanging $r$ and $s$ is also satisfied, and the identity matrix $1_n$ is in the image of $A(i,i)$ for all $i = 1,\dots,m$, and the dual condition exchanging $i$ and $j$ is also satisfied.

There are three main families of twisting maps arising in this way. The \emph{standard} twisting maps $\tau$ have only entries $0$ or $1$ in all the columns of $M_{\tau}(s,s)$, and the non-zero entries of $M_{\tau}(i,s)$ are at most the non-zero entries of $M_{\tau}(s,s)$, for all $i,s = 1,\dots, m$. These give algebras isomorphic to certain square zero radical truncated quiver algebras. The \emph{quasi-standard} twisting maps $\tau$ arise as generalizations of the standard ones where some of the conditions are relaxed, and give algebras that are formal deformations of the standard ones. The third family is given only when $n = m$, whenever there are $v_1,\dots v_n \in \kk^n$ invertible such that $\det(v_1^T,\dots,v_n^T) = 1$, in which case for all $i,s = 1,\dots,n$
\begin{equation*}
M_{\tau}(i,s) = (-1)^{i-1}(v_s^{-1} \cdot v_i)^T (v_l \cdot (v_1 \times \cdots \times \widehat{v_i} \times \cdots \times v_n))
\end{equation*}
and these give algebras isomorphic to the square matrices over the field $\kk$. In general these three families of twisting maps do not give a complete description of all the possibilities, but for $\kk^3$ with $\kk^3$ the classification of the remaining ones was completed on a case-by-case basis \cite{MR4283170}.

The study of twisted planes exploits Koszulness and Artin-Schelter regularity, and techniques of linear algebra involving infinite matrices give an almost complete classification \cite{BancesValqui}.

\begin{theorem} \cite[Theorem 1.3]{MR4248207}
Let $\kk$ be an algebraically closed field. A quadratic twisted tensor product $\kk[x]\otimes_{\tau} \kk[y]$ is isomorphic as a graded twisted tensor product to
\begin{equation*}
\kk\langle x,y\rangle/(yx - ax^2 - bxy - y^2)\quad \text{or}\quad \kk\langle x,y\rangle/(yx - ax^2 - bxy)
\end{equation*}
for some $a,b\in \kk$.
\end{theorem}

Note that the notion of isomorphism as graded twisted tensor products is stronger than that of isomorphism as graded algebras.

In general, little more is classified besides the above. A second typical approach would be to inquire whether twisted tensor products inherit any properties from the given algebras. For example, since connected Ore extensions preserve Artin-Schelter regularity \cite[Theorem 0.2]{MR2452318}, one may ask what other properties are preserved. In general, this is a hard question. There is an appropriate sense in which Artin-Schelter regularity is preserved \cite{MR3866681}.

\begin{example} \label{ex.ttp}
Let $A = \kk[u, v]$ and $B = \kk[x, y]$ with $|u| = |v| = |x| = |y| = 1$. Consider the strongly graded twisting map $\tau: B \otimes A \to A \otimes B$ given by 
\[\tau(x \otimes u) = a u \otimes x,\quad \tau(x \otimes v) = b v \otimes x,\quad \tau(y \otimes u) = c u \otimes y,\quad \tau(y \otimes v) = d v \otimes y\]
for some $a,b,c,d\in \kk^{\times}$. The twisted tensor product $A \otimes_{\tau} B$ is the free algebra on generators $u, v, x$, and $y$ subject to the following relations
\begin{align*}
&uv - vu,  && aux - xu, && cuy - yu,\\
&bvx - xv, && dvy - yv, && xy - yx,
\end{align*}
where we slightly abuse notation and write $u=u\otimes 1, v=v\otimes 1, x=1\otimes x, y=1\otimes y$. Observe that $A\otimes_{\tau} B$ has a PBW-type decomposition, namely it can be identified with the $\kk$-vector space having basis $\{u^i v^j x^r y^s\}_{i,j,r,s\in\mathbb{N}}$ and multiplication
\begin{equation*}
u^{i_1} v^{j_1} x^{r_1} y^{s_1} \cdot u^{i_2} v^{j_2} x^{r_2} y^{s_2} = a^{i_1 i_2} b^{j_1 j_2} c^{r_1 r_2} d^{s_1 s_2} u^{i_1+i_2} v^{j_1+j_2} x^{r_1+r_2} y^{s_1+s_2}
\end{equation*}
for all $i_1,i_2,j_1,j_2,r_1,r_2,s_1,s_2\in\mathbb{N}$. In particular, $A\otimes_{\tau} B$ inherits a $\Z$-grading by setting $|u^i v^j x^r y^s| = i+j+r+s$ and $A \otimes_{\tau} B$ is a quadratic Artin-Schelter regular algebra of global dimension four. If $a = b = c = d$, then $A \otimes_{\tau} B$ is given by a Zhang twist of the polynomial ring $\kk[u,v,x,y]$. Indeed, taking the graded algebra automorphism $\theta$ of $\kk[u,v,x,y]$ defined by
\[\theta(u) = au,\quad \theta(v) = av,\quad \theta(x) = x,\quad \theta(y) = y,\]
yields $A \otimes_{\tau} B \cong \kk[u,v,x,y]^{\theta}$ (see Example \ref{exam.zhang-a}). Conversely, if $A \otimes_{\tau} B$ is a Zhang twist of $\kk[u, v,x,y]$, then $bd^{-1} = ac^{-1} = dc^{-1} = ba^{-1}=1$ by \cite[Lemma 2.2]{NV}, that is $a = b = c = d$. Therefore, Zhang twists and twisted tensor products are different but related kinds of twists, and have been shown to be useful in the study of Artin-Schelter regular algebras.
\end{example}

It is also known that twisted tensor products of Hopf algebras never inherit a Hopf algebra structure (as long as the twisting map is not trivial), but that twisted tensor products of Frobenius algebras always inherit a Frobenius algebra structure \cite{OcalOswald}. However, the following example shows that the twisted tensor product of Koszul algebras need not be Koszul.

\begin{example} \cite[Example 5.4]{MR3849878}
Consider two polynomial rings $\kk[x]$ and $\kk[y]$. We can define a $\kk$-linear map $\tau:\kk[y]\otimes \kk[x] \to \kk[x]\otimes \kk[y]$ by setting $\tau(y\otimes x) = q x\otimes y$ for some non-zero $q\in \kk$, and extending as before. This gives a strongly graded twisting map by construction, and the twisted tensor product is isomorphic to the quantum plane $\kk_q[x,y] = \kk\langle x,y\rangle/(qxy + yx)$, which is Koszul.

If instead we consider $\tau:\kk[y]\otimes \kk[x] \to \kk[x]\otimes \kk[y]$ by setting $\tau(y^i\otimes x^j) = x^j\otimes y^i$ when $i$ or $j$ are even, $\tau(y^i\otimes x^j) = x^{j+1}\otimes y^{i-1} - x^j \otimes y^i + x^{j-1} \otimes y^{i+1}$ when $i$ and $j$ are both odd, then the twisted tensor product is not quadratic (whence not Koszul) and it has infinite global dimension (whence it is not Artin-Schelter regular).
\end{example}

Generalizations of the former construction where the twist only involves multiplying the graded generators by a scalar have been extensively considered. They are known as \emph{twists by a bicharacter}, the resulting twisted tensor product is denoted by $A\otimes^t B$. These are constructed from algebras $A = \bigoplus_{f\in F}{A_f}$ and $B = \bigoplus_{g\in G}{B_g}$ graded by abelian groups $F$ and $G$. Let $t:F\otimes_{\mathbb{Z}} G \to \kk^{\times}$ be a $\kk$-bicharacter, namely
\begin{align*}
t(1_F,1_G) = 1,\quad t(f_1 f_2,g) = t(f_1,g)t(f_2,g),\quad t(f,g_1g_2) = t(f,g_1)t(f,g_2),
\end{align*}
then $\tau:B\otimes A\to A\otimes B$ given by linearly extending $\tau(b\otimes a) = t(|a|,|b|) a\otimes b$ for all homogeneous $a\in A$ and $b\in B$ is a twisting map. The main example are quantum complete intersections
\begin{equation*}
\kk\langle x_1,\dots, x_n\rangle/\left(x_i^{m_i},x_i x_j - q_{ij} x_j x_i \right)_{i,j\in\{1,\dots,n\}}
\end{equation*}
where $n, m_1,\dots,m_n\in\N$, $n,m_1,\dots,m_n \geq 2$, $(q_{ij}) \in M_{n}(\kk^{\times})$, $q_{ii} = 1$, and $q_{ij} q_{ji} = 1$, who have been extensively studied for their relation to Happel's question \cite{MR2189240}, interesting support theories \cite{MR2341266, MR2563183}, and (co)homologies \cite{MR2429451, MR3879091, MR2776874}. In fact, twists by a bicharacter were introduced in \cite{MR2450729} with a view to generalize the fact that the Hochschild cohomology of a tensor product is the tensor product of the Hochschild cohomologies, a result formalized in \cite{BriggsWitherspoon}.

\begin{theorem} \cite[Theorem 3.3]{MR3175033}
Let $A$ and $B$ be algebras, one of them finite dimensional. Then there is an isomorphism of Gerstenhaber algebras
\begin{equation*}
\HH^{\bullet}(A\otimes B) \cong \HH^{\bullet}(A)\otimes \HH^{\bullet}(B).
\end{equation*}
\end{theorem}

\begin{theorem} \cite[Theorem 5.1]{BriggsWitherspoon}
Let $A$ and $B$ be algebras graded by abelian groups $F$ and $G$. Then there is an isomorphism of Gerstenhaber algebras
\begin{equation*}
\HH^{\bullet}(A\otimes^t B) \cong \bigoplus_{f\in F,\, g\in G}{\HH^{\bullet}(A,A_{\hat{g}})_f\otimes \HH^{\bullet}(B,\prescript{}{\hat{f}}B)_g}
\end{equation*}
where $\hat{f}$ and $\hat{g}$ are automorphisms of $B$ and $A$ given by the bicharacter, $A_{\hat{g}}$ is the graded $A$ bimodule with twisted right multiplication, and $\prescript{}{\hat{f}}B$ is the graded $B$ bimodule with twisted left multiplication.
\end{theorem}
There is a method to compute the Hochschild cohomology of a general twisted tensor product in terms of its components. Fix algebras $A$, $B$, and a twisting map $\tau$. Given $(M,\nabla_M)$ an $A$ bimodule and $(N,\nabla_N)$ a $B$ bimodule together with bijective $\kk$-linear maps $\tau_{M} : B\otimes M\to M\otimes B$ and $\tau_{N} : N\otimes A\to A\otimes N$ such that
\begin{align*}
\tau_{M} (\nabla_B \otimes 1) &= (1\otimes \nabla_B)(\tau_M\otimes 1) (1\otimes \tau_M),\\
\tau_{N} (1\otimes \nabla_A) &= (\nabla_A \otimes 1)(1\otimes \tau_N) (\tau_N \otimes 1),\\
\tau_{M} (1 \otimes \nabla_M) &= (\nabla_M\otimes 1)(1\otimes 1\otimes \tau) (1\otimes \tau_M \otimes 1) (\tau \otimes 1\otimes 1),\\
\tau_{N} (\nabla_N \otimes 1) &= (1\otimes \nabla_N)(\tau\otimes 1\otimes 1) (1\otimes \tau_N \otimes 1) (1\otimes 1 \otimes \tau),
\end{align*}
then $M\otimes N$ is an $A\otimes_{\tau} B$ bimodule, denoted $(M\otimes_{\tau} N,\nabla_{M\otimes_{\tau} N})$, where
\begin{equation*}
\nabla_{M\otimes_{\tau} N} \coloneqq (\nabla_M \otimes \nabla_N) (1\otimes 1\otimes \tau \otimes 1\otimes 1) (1\otimes \tau_{M}\otimes \tau_{N}\otimes 1).
\end{equation*}
In this case, $M$ and $N$ are said to be \emph{compatible} with $\tau$. Given resolutions of $M$ and $N$, they are said to be \emph{compatible} with $\tau$ when each of their modules is compatible, and the collection of resulting $\kk$-linear maps lift $\tau_M$ and $\tau_N$.

\begin{theorem} \cite[Theorem 3.10]{MR3936025}
Let $A$ and $B$ be algebras, let $\tau:B\otimes A\to A\otimes B$ be a twisting map, let $M$ and $N$ be $A$ and $B$ bimodules respectively, let $P_{\bullet}$ and $Q_{\bullet}$ be free resolutions of $M$ and $N$ as $A$ and $B$ bimodules respectively. Assume that $M$, $N$, $P_{\bullet}$, and $Q_{\bullet}$ are all compatible with $\tau$. Then the tensor product complex $(P_{\bullet} \otimes Q_{\bullet})_{\bullet}$ is a free resolution of $M\otimes N$ as a $A\otimes_{\tau} B$ bimodule.
\end{theorem}

The above, together with the fact that the bar resolution is compatible with all twisting maps, allows the computation of $\HH^{\bullet}(A\otimes_{\tau} B)$ in general. The Gerstenhaber algebra structure can be found using the tools introduced in \cite{MR4168969}.

\section{Twisted Segre Products}
\label{Sec4}

From now on we assume for simplicity that $\kk$ is algebraically closed. The Segre product of graded algebras plays an important role in algebraic geometry and commutative algebra (see e.g., \cite{Harris, GW}). In this section, we discuss the notion of twisted Segre product of graded algebras.

\begin{definition}
Let $A=\bigoplus_{i \in \mathbb Z}A_i$ and $B=\bigoplus_{i \in \mathbb Z}B_i$ be (not necessarily commutative) $\mathbb Z$-graded algebras. The \emph{Segre product} $A\circ B$ is the $\kk$-vector space
\begin{align*}
A \circ B \coloneqq \bigoplus_{i \in \mathbb Z} (A_i \otimes B_i)
\end{align*}
with the multiplication map
\[
\nabla_{A\circ B} := (\nabla_A\otimes \nabla_B)(1\otimes \sigma \otimes 1),
\]
where $\sigma : B\otimes A \to A\otimes B$ transposes the entries. This makes $A\circ B$ a $\Z$-graded algebra.
\end{definition}

Let $X \subseteq {\mathbb P}^{n-1}$ and $Y \subseteq {\mathbb P}^{m-1}$ be projective varieties with the homogeneous coordinate rings $A$ and $B$, respectively. Consider the projective variety $X \times Y$ in ${\mathbb P}^{nm-1}$ via the Segre embedding. Then $A \circ B$ is the homogeneous coordinate ring of $X \times Y$.

\begin{definition} \cite[Section 3]{HU}
Let $A, B$ be $\mathbb Z$-graded algebras and let $\tau: B\otimes A\to A\otimes B$ be a strongly graded twisting map. The \emph{twisted Segre product} $A\circ_\tau B$ is the $\kk$-vector space $A\circ B$ with the multiplication map
\[
\nabla_{A\circ_{\tau} B} \coloneqq (\nabla_A\otimes \nabla_B)(1\otimes \tau \otimes 1).
\]
This makes $A\circ_{\tau} B$ a $\Z$-graded algebra.
\end{definition}

We give some further definitions. Let $S$ be a Noetherian $\mathbb Z$-graded algebra, denote by $\operatorname{GrMod} S$ the category of graded right $S$-modules. We say that an element $m \in M \in \operatorname{GrMod} S$ is \emph{torsion} if $\dim_\kk mS< \infty$. A graded module $M \in \operatorname{GrMod} S$ is called \emph{torsion} if every homogeneous element of $M$ is torsion.
We denote by $\operatorname{Tors} S$ the full subcategory of $\operatorname{GrMod} S$ consisting of torsion modules and by $\operatorname{QGr} S$ the Serre quotient category $\operatorname{GrMod} S/\operatorname{Tors} S$. The category $\operatorname{QGr} S$ is frequently thought of as the category of quasi-coherent sheaves on the imaginary projective scheme associated to $S$ (see \cite{AZ}).

To study $\operatorname{QGr} (A\circ_{\tau} B$), where $A$ and $B$ are connected graded algebras finitely generated in degree $1$, we use bigraded algebras. Let $S=\bigoplus_{i,j \in \mathbb Z}S_{(i,j)}$ be a $\mathbb Z \times \mathbb Z$-bigraded algebra, define $S_i \coloneqq \bigoplus_{j \in \mathbb Z}S_{(i,j)}$ for each $i \in \mathbb Z$. Then $S=\bigoplus_{i \in \mathbb Z}S_{i}$ can be viewed as a $\mathbb Z$-graded algebra and $S_0= \bigoplus_{j \in \mathbb Z}S_{(0,j)}$ is a $\mathbb Z$-graded algebra.

\begin{definition} \cite[Section 4.1]{HU}
\begin{enumerate}
    \item A $\mathbb Z$-graded algebra $S=\bigoplus_{i \in \mathbb Z}S_i$ is called a \emph{densely graded algebra} when $\dim_\kk (S_{i+j}/S_{i}S_{j})<\infty$ for all $i, j \in \mathbb Z$.
    \item A $\mathbb Z \times \mathbb Z$-bigraded algebra $S=\bigoplus_{i,j \in \mathbb Z}S_{(i,j)}$ is called a \emph{densely bigraded algebra} when the $\mathbb Z$-graded algebra $S=\bigoplus_{i \in \mathbb Z}S_{i}$ is a densely graded algebra, where $S_i = \bigoplus_{j \in \mathbb Z}S_{(i,j)}$.
\end{enumerate}
\end{definition}
A $\mathbb Z$-graded algebra $S$ is said to be \emph{strongly graded} when $S_{i}S_{j}=S_{i+j}$ for all $i, j \in \mathbb Z$. Strongly graded algebras are clearly densely graded algebras. Moreover if $S$ is a strongly graded algebra then the functor $\operatorname{GrMod} S \to \operatorname{Mod} S_0$ defined by $M \mapsto M_0$ is an equivalence by Dade's Theorem \cite[Theorem 3.1.1]{NV2004}.

Let $S$ be a densely bigraded algebra satisfying the following conditions.
\begin{enumerate}
\item[{(D1)}] $S_0\coloneqq \bigoplus_{j \in \mathbb Z}S_{(0,\,j)}$ is a Noetherian $\mathbb Z$-graded algebra, and
\item[{(D2)}] $S_i\coloneqq \bigoplus_{j \in \mathbb Z}S_{(i,\,j)}$ is a finitely generated graded right and left $S_0$-module for every $i \in\mathbb Z$.
\end{enumerate}
Denote by $\operatorname{BiGrMod}S$ the category of bigraded right $S$-modules. We say that an element $m \in M \in\operatorname{BiGrMod} S$ is \emph{locally torsion} if $\dim_\kk mS_0< \infty$.
A bigraded module $M \in \operatorname{BiGrMod} S$ is called \emph{locally torsion} if every homogeneous element of $M$ is locally torsion.
We denote by $\operatorname{BiLTors} S$ the full subcategory of $\operatorname{BiGrMod} S$ consisting of locally torsion modules and by
$\operatorname{QBiGr_{L}} S$ the Serre quotient category $\operatorname{BiGrMod} S/ \operatorname{BiLTors} S$.
We then have the following version of Dade’s Theorem for densely bigraded algebras.

\begin{theorem} \cite[Theorem 4.10]{HU} \label{thm.DT}
Let $S$ be a densely bigraded algebra satisfying conditions \textnormal{(D1)} and \textnormal{(D2)}.
Then the functor $(-)_0: \operatorname{BiGrMod} S \to \operatorname{GrMod} S_0$ given by $M \mapsto M_0 \coloneqq \bigoplus_{j \in \mathbb Z}M_{(0,j)}$ induces an equivalence
\[ \operatorname{QBiGr_{L}} S \cong \operatorname{QGr} S_0. \]
\end{theorem}

Let us go back to consider twisted Segre products. Let $A$ and $B$ be $\mathbb Z$-graded algebras and let $\tau: B\otimes A\to A\otimes B$ be a graded twisting map. We endow the twisted tensor product $A \otimes_\tau B$ with a bigraded structure by setting
\[ (A \otimes_\tau B)_{(i,j)} \coloneqq A_{i+j}\otimes B_{j}\]
for each $i,j \in \mathbb Z$.
Then $(A \otimes_\tau B)_0$ is equal to the twisted Segre product $A\circ_\tau B$ as graded algebras.

\begin{proposition} \cite[Proposition 4.11]{HU} \label{prop.dbga}
Let $A$ and $B$ be connected graded algebras finitely generated in degree $1$ and let $\tau: B \otimes A \to A\otimes B$ be a strongly graded twisting map. Then $A\otimes_\tau B$ is a densely bigraded algebra satisfying condition \textnormal{(D2)}.
\end{proposition}

Theorem \ref{thm.DT} and Proposition \ref{prop.dbga} yield the following result which gives a categorical relationship between twisted tensor products and twisted Segre products.

\begin{theorem} \cite[Theorem 4.13]{HU}\label{Segremainthm}
Let $A$ and $B$ be connected graded algebras finitely generated in degree $1$ and let $\tau: B \otimes A \to A\otimes B$ be a strongly graded twisting map. If $A \circ_\tau  B$ is Noetherian then there is an equivalence
\[ \operatorname{QBiGr_{L}} (A \otimes_\tau  B) \cong \operatorname{QGr} (A \circ_{\tau} B). \]
\end{theorem}

\begin{remark}
\begin{enumerate}
    \item In general, it is not known whether the twisted Segre product $A \circ_\tau  B$ is Noetherian when $A$ and $B$ are Noetherian.
    As a positive example, the usual Segre product $A \circ B$ is Noetherian if $A$ and $B$ are commutative finitely generated connected graded algebras.
    \item A category equivalence similar to that of Theorem \ref{Segremainthm} was proved in \cite[Theorem 2.4]{VR} for the usual Segre product. However, the bigradings of the tensor products $A\otimes B$ given here differ from the ones in \cite{VR}, and hence Theorem \ref{Segremainthm} is not a generalization of \cite[Theorem 2.4]{VR}.
\end{enumerate}

\end{remark}

If $A$ and $B$ are Noetherian Koszul Artin-Schelter regular algebras and $\tau: B \otimes A \to A\otimes B$ is a strongly graded twisting map, then it follows from \cite[Theorem 2]{MR3866681} that $A \otimes_\tau  B$ has finite global dimension. From this we see that $\operatorname{QBiGr_{L}} (A \otimes_\tau  B)$ has finite global dimension. Therefore, as a corollary of Theorem \ref{Segremainthm}, we have the following.

\begin{corollary} \cite[Theorem 4.16]{HU}\label{Segremaincor}
Let $A$ and $B$ be Noetherian Koszul Artin-Schelter regular algebras and $\tau: B \otimes A \to A\otimes B$ be a strongly graded twisting map. If $A \circ_\tau  B$ is Noetherian then $\operatorname{QGr} (A \circ_{\tau} B)$ has finite global dimension.
\end{corollary}

The above corollary can be considered as a noncommutative analogue of the fact that a Segre variety $\mathbb P^{n-1} \times \mathbb P^{m-1} \cong \operatorname{Proj}(\kk[x_1,\dots,x_n]\circ \kk[y_1,\dots,y_m])$ is smooth.

\begin{example}
We now continue Example \ref{ex.ttp}. Consider the (ungraded) subalgebra of $A\otimes_{\tau} B$ generated by the basis elements $ux, vx, uy$, and $vy$, where we write $ux=u\otimes x$, $vx=v\otimes x$, $ uy=u\otimes y$, and $vy=v\otimes y$, slightly abusing notation as before. Their multiplication is given by
\begin{align*}
&a ux\cdot vx = b vx\cdot ux,  && c ux\cdot uy = a uy\cdot ux, && c ux\cdot vy = b vy\cdot ux,\\
&d vx \cdot uy = a uy \cdot vx, && d vx \cdot vy = b vy\cdot vx, && c uy\cdot vy = d vy\cdot uy,\\
& a ux\cdot vy = b vx\cdot uy. && &&
\end{align*}
Denote $X = uv$, $Y= vx$, $Z= uy$, $W= vy$ and regrade $|X|=|Y|=|Z|=|W|=1$. The twisted Segre product $A \circ_{\tau} B$ is the $\mathbb{Z}$-graded algebra freely generated by $X, Y, Z$, and $W$ subject to the relations
\begin{align*}
&aXY-bYX, &&cXZ-aZX, &&cXW-bWX,\\
&dYZ-aZY, &&dYW-bWY, &&cZW-dWZ, &&aXW-bYZ.
\end{align*}
It is straightforward to check that $A \circ_{\tau} B$ is a quotient ring of a skew polynomial ring, so it is Noetherian. Hence $\operatorname{QGr} (A \circ_{\tau} B)$ has finite global dimension by Corollary \ref{Segremaincor}.
\end{example}

An interesting future research topic would be to further explore properties of twisted Segre products.

\section*{Acknowledgments} We thank the organizers of the conferences ``Recent Advances and New Directions in the Interplay of Noncommutative Algebra and Geometry" (Seattle 2022) and ``Noncommutative Geometry and Noncommutative Invariant Theory'' (BIRS 2022) for enabling this collaboration. In particular, we effusively thank Michaela Vancliff and James J. Zhang for their continued support. The first author was supported by an AMS-Simons travel grant, and by the Hausdorff Research Institute for Mathematics funded by the Deutsche Forschungsgemeinschaft (DFG, German Research Foundation) under Germany's Excellence Strategy - EXC-2047/1 - 390685813. The second author was supported by JSPS KAKENHI Grant Numbers JP18K13381 and JP22K03222.

\bibliography{bibTSurvey}
\bibliographystyle{plain}
\end{document}